\documentclass[12pt]{amsart}

\usepackage[T1]{fontenc}
\usepackage{comment}
\usepackage{appendix}
\usepackage{graphicx}
\usepackage{pdfpages}
\newtheorem{theorem}{Theorem}[section]
\theoremstyle{definition}
\newtheorem{definition}[theorem]{Definition}
\theoremstyle{remark}
\newtheorem{remark}[theorem]{Remark}
\usepackage{subcaption}

\usepackage{epic}
\usepackage{import}
\usepackage{pdfpages}

\usepackage{float}

\title{New Exotic four-manifolds with $\mathbb{Z}/2\mathbb{Z}$ fundamental group}
\author{M\'{a}rton Beke, L\'{a}szl\'{o} Koltai, and Sarah Zampa}

\newcommand{\ie}{\textit{i}.\textit{e}., }
\def\Z{\mathbb{Z}}
\def\Q{\mathbb{Q}}
\DeclareMathOperator{\CP}{\mathbb{CP}}

\begin{document}

\begin{abstract}
We extend a construction of Stipsicz-Szab\'{o}  of infinitely many irreducible exotic smooth structures of some closed four-manifolds with even $b_2^+$ and fundamental group $\mathbb{Z}/2\mathbb{Z}$.
We use the double node surgery and rational blow down constructions of Fintushel-Stern on some elliptic fibrations equipped with a free involution.
The construction is done in an equivariant manner and the factor manifolds are distinguished by the Seiberg-Witten invariants of their universal covers.
\end{abstract}
{{\maketitle}}

\section{Introduction}

Two smooth manifolds are said to be \textit{exotic} if they are homeomorphic but not diffeomorphic to each other. In this case, we call the two smooth structures \textit{exotic}, and a manifold invariant is usually used to distinguish smooth structures, \ie to show non-diffeomorphism of the two manifolds.
In the following, we use the Seiberg-Witten invariants \cite{witten1994monopoles} given by the function
\[
SW_X : H^2(X;\Z) \rightarrow \mathbb{Z}
\]
where $X$ is a smooth, closed, oriented four-manifold with $b_2^+>1$.
This is encoded in a formal series $\mathcal{SW}_X:=\sum_{\alpha\in H^2(X;\Z)}SW_X(\alpha)e^\alpha$.
Note, that the number of \textit{``basic classes''}, \textit{i.e.} classes $\alpha$ with $SW_X(\alpha)\neq 0$, is finite, and so the sum is also finite.
Indeed, if two manifolds are diffeomorphic, then their Seiberg-Witten invariants have to agree (up to automorphism).
On the other hand, the homeomorphism of the two manifolds will arise from the generalization of Freedman's classification theorem to manifolds with fundamental group $\mathbb{Z}/2\mathbb{Z}$ (see \cite{hambleton1993cancellation}).

In this paper, we will expand the construction given in \cite[Theorem 4.6]{stipsicz2023exotic} by decomposing the elliptic fibration
\[E(2n+1)=E(n)\#_fE(1)\#_fE(n)\]
and exploiting the $E(n)$ summands of the fibration.
We use the rational blow down construction, knot and double node surgery constructions of Fintushel-Stern (see \cite{fintushel1997rational, fintushel1998knots, fintushel2006double} respectively) in order to construct infinite families of exotic simply connected four manifolds. Those manifolds will come equipped with involutions, which will provide the exotica of the following theorem:

\begin{theorem}
    The manifold $Z_1\#2n\CP^2\#l\overline{\CP^2}$ where $l\in\{5n+6,5n+12,\dots,8n\}$ if $n=2q\not\equiv6\pmod{8}$ and $l\in\{5n+9,5n+15,\dots,8n\}$ if $1<n=2q+1\not\equiv 5\pmod{8}$, admits infinitely many {irreducible} smooth structures.
\end{theorem}

Here, $Z_1$ denotes the quotient of $S^2\times S^2$ with the fixed point free involution $\iota$, which applies the antipodal map on the two $S^2$ components.
We conjecture that the same theorem is valid in the remaining $n\equiv 5,6\pmod{8}$ cases as well.

\vspace{5mm}
\noindent\textbf{Acknowledgements:} We wish to thank Andr\'{a}s Stipsicz for his help and guidance throughout this project.

\section{Preliminaries}

\subsection{Topological results}
We will make use of the homeomorphism classification of oriented smooth manifolds with $\pi_1=\mathbb{Z}/2\mathbb{Z}$, to this end consider $X$ and its two-fold universal cover $\tilde X$.
There are 3 possibilities with regards to the spinness of these two manifolds:
\begin{enumerate}
    \item[I)] neither are spin
    \item[II)] both are spin
    \item[III)] $X$ is not spin, but $\tilde X$ is.
\end{enumerate}
This is called the \textit{$w_2$-type} of the manifold $X$.
Note, that these are the only possibilities, since spinness is equivalent with having $w_2=0$, and naturality of this class means, that if $X$ is spin, then so is $\tilde X$.
Now we can state 
\begin{theorem}[{\cite[Theorem C]{hambleton1993cancellation}}]
    Let $X_1,\ X_2$ be two closed oriented smooth 4-manifolds with fundamental group $\mathbb{Z}/2\mathbb{Z}$.
    The manifolds are homeomorphic if and only if their Euler characteristic, signature and $w_2$-type agree.
\end{theorem}

\subsection{Knot surgery}
Consider a knot $K\subset S^3$, and a simply connected smooth manifold $X$ with $b_2^+>1$, and with an embedded homologically essential torus $T$ 
of self-intersection $0$ and simply connected complement.
The manifold $X_K:=(X\setminus \nu T)\cup_\phi(S^3\setminus\nu K)\times S^1$ is called a knot surgery of $X$ using $K$, where $\phi:T^3\to T^3$ is chosen so that the longitude of $K$ is identified with the normal circle of $T$, and $T^3=S^1\times S^1\times S^1$ denotes the $3$-torus. Note that this does not determine $\phi$ completely: we pick and fix such a function.

The following theorem shows the importance of this construction:
\begin{theorem}[\cite{fintushel1998knots}]
    With the setup as above $\mathcal{SW}_{X_K}=\mathcal{SW}_X\cdot\Delta_K(e^{2[T]})$ holds, where $\Delta_K$ is the symmetric Alexander polynomial of $K$.
\end{theorem}

In the special case of elliptic fibrations over the sphere -- which we will be using -- more can be done (\cite{fintushel2006double}). Let $K$ be any genus one knot, and pick a minimal genus Seifert surface $\Sigma$  
and a 
non-separating loop $\Gamma$ inside $\Sigma$ satisfying:
\begin{enumerate}
    \item $\Gamma$ bounds a disc which intersects $K$ in two points.
    \item $\Gamma$ has linking number $+1$ in $S^3$ with its pushoff on $\Sigma$.
\end{enumerate}
\begin{definition}
Consider a neighborhood of a smooth fiber $F$ in an elliptic fibration containing exactly two nodal fibers with the same monodromy, and with all other fibers smooth.
Such a neighborhood is called a \textit{double node neighborhood}.
\end{definition}

With $K$ as previously and $T$ a smooth fiber in a double node neighborhood, we do knot surgery with the map $\phi$ picked so that the meridian of $K$ is identified with the vanishing cycle of the 
nodal fibers (both are the same, since they have the same monodromy).

In the local picture a section of the elliptic fibration looks like a disc; 
the result of knot surgery is to remove a smaller disc from the section  
and replace it with the Seifert surface of $K$, 
thus, we obtain a punctured torus afterwards. Inside this torus sits the loop $\Gamma$, which by $(1)$ bounds a twice punctured disc $D'$ with boundary $\partial D'=\Gamma \cup m_1 \cup m_2$, where the $m_i$ are meridians of $K$.
By the choice of $\phi$ the meridians $m_i$ bound disjoint vanishing discs $D_i$ (corresponding to the two nodal singularities), hence $\Gamma$ bounds a disc $D=D'\cup D_1 \cup D_2$. With respect to the framing of $\Gamma$ given by the push-off on the surface $\Sigma$, the relative self-intersection of $D$ is $-1$.

Remove now an annular neighborhood of $\Gamma$ inside $\Sigma$, and close the resulting surface with two copies $B_1$ and $B_2$ of the disc $D$ described above. These two discs intersect each other in one point, so the capped off surface is an immersed disc with a double point. To get the sign of the double point notice that $B_1$ and $B_2$ intersect with negative sign, but we have to change the orientation of $B_2$ in order to get an oriented surface, so the immersed disc has a positive double point. Thus in a double node neighborhood we can replace the genus introduced by the knot surgery with a positive double point if we use knots satisfying the conditions above.

\subsection{Rational blow down}
Another ingredient is the rational blow down (\cite{fintushel1997rational,park1997seiberg}).
Denote by $C_p$ a linear plumbing (consecutive terms intersect transversely at a single point) of $p-1$ spheres with self-intersections $-(p+2),-2,\dots,-2$.
The boundary of this linear plumbing is the lens space $\partial C_p=L(p^2,p-1)$.
The main observation is that there is a manifold $B_p$ with the same boundary, and with $H_*(B_p;\Q)=H_*(\{pt\};\Q)$ (\cite[Figure 8.42.]{gompf19994}).
If $C_p\subset X$ for some simply connected 4-manifold $X$, we call $X_p:=(X\setminus\nu C_p)\cup B_p$ the rational blow down of $X$.

\begin{figure}[ht!]
    \centering
\begin{center}
\begin{picture}(30,20)(90,60)
\put(50,60){\circle*{5}}
\put(90,60){\circle*{5}}
\put(130,60){\circle*{5}}
\put(170,60){\circle*{5}}
\put(50,60){\line(1,0){80}}

\put(45,70){\makebox(0,0){\small{$-(p+2)$}}}
\put(90,70){\makebox(0,0){\small{$-2$}}}
\put(130,70){\makebox(0,0){\small{$-2$}}}
\put(170,70){\makebox(0,0){\small{$-2$}}}
\dashline{2}(130,60)(170,60)
\put(50,50){\makebox(0,0){\small{$u_1$}}}
\put(90,50){\makebox(0,0){\small{$u_2$}}}
\put(130,50){\makebox(0,0){\small{$u_3$}}}
\put(170,50){\makebox(0,0){\small{$u_{p-1}$}}}

\end{picture}
\end{center}

    \caption{The configuration $C_p$}
    \label{fig:Cp}
\end{figure}
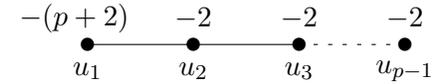

In special circumstances it is possible to follow the Seiberg-Witten invariants after a rational blow down.
It is true in general, that if $X_p$ is simply connected, then the Seiberg-Witten invariants of $X$ completely determine those of $X_p$, but the computation is not always simple.
A more manageable case which we will need is the following:

\begin{definition}\label{def:taut}
    $C_p$ is \textit{Seiberg-Witten tautly embedded} if $|\alpha(u_1)|\leq p$ (where $u_1\subset C_p$ is the sphere with self-intersection $-p-2$), and $\alpha(u_i)=0$ (for the ($-2$)-spheres of the configuration) is satisfied for all basic classes $\alpha$.
\end{definition}

To obtain the Seiberg-Witten invariants after the blow down, one needs to find extensions of certain cohomology classes from $X\setminus\nu C_p$ to $X_p$, and vice versa.
If an extension exists, and the expected dimension of the moduli space (defined by the right hand side of Equation~\eqref{eq:smiple}) stays non-negative, then the value of the invariant is unchanged after the blow down \cite[Theorem 8.2]{fintushel1997rational}. In our case, the expected dimension of the moduli space is always zero.

In the tautly embedded case this extension process simplifies further, and one only has to consider the basic classes which evaluate maximally on $u_1$:
\begin{theorem}
    Let $X$ be a simply connected 4-manifold with $b_2^+\geq 2$, and $C_p$ tautly embedded into it.
    If $\alpha'$ is a basic class of $X_p$, and $\alpha|_{X\setminus\nu C_p}=\alpha'|_{X_p\setminus B_p}$, then $|\alpha(u_1)|=p$ where $u_1\subset C_p$ is the sphere with $[u_1]^2=-p-2$.
\end{theorem}
Finally the standard blow up formula also needs to be mentioned, as we will be relying on it in the following.
\begin{theorem}[{\cite[Theorem 1.4]{fintushel1995immersed}}]
    Let $X$ be a 4-manifold of Seiberg-Witten simple type with $b_2^+>1$, and $X'=X\#\overline{\CP^2}$ its blow up, then $SW_X(\alpha)=SW_{X'}(\alpha\pm E)$, with the (Poincaré dual of the) new exceptional sphere denoted $E\in H^2(X\#\overline{\CP^2};\Z)$.
\end{theorem}
\begin{remark}
    For the reader unfamiliar with Seiberg-Witten invariants, the \textit{simple type} condition, i.e. that for every basic class $\alpha$ \begin{equation}\label{eq:smiple}
        \frac{{\alpha^2([X])}-(3\sigma(X)+2\chi(X))}{4}=0
    \end{equation} can be safely ignored, since $E(n)$ is of simple type for all $n>1$ and simple type manifolds stay simple type under knot surgeries and rational blow downs.
    For an exposition, see e.g. \cite[Section 8]{fintushel1997rational}.
\end{remark}

\subsection{Monodromy factorisation.} \label{monodromy}
Let $f:X\rightarrow Y$ be an elliptic fibration, and $\Delta=\{p_1,\dots,p_n\}\subset Y$ be the discriminant locus of the fibration (i.e. the points whose fiber is singular). 
We assume that $Y$ is the sphere $S^2$. 
Let $p\in \Delta$, and restrict the fibration to a small circle around $p$, which is disjoint from $\Delta$.
By going around the small circle, we get a self-diffeomorphism of the generic fiber $T^2$, defined up to isotopy and conjugation.
In other words, we get an element defined up to conjugation in the mapping class group $\Gamma_1$ of $T^2$.
We will call this conjugacy class the \textit{monodromy} of the singular fiber $f^{-1}(p)$.
The mapping class group of the torus admits a presentation:
 $$\Gamma_1=\langle a,b \ | \ aba=bab, (ab)^6=1 \rangle, $$
 where $a$ and $b$ correspond to the right-handed Dehn twists around the natural generators of $H_1(T^2)$.
 The monodromy of the singular fiber $I_n$ is conjugate to $a^n$ in this presentation.
 Throughout this paper we use the convention $a^b=\overline{b}ab$ where $a\overline a=1$.

 Suppose now that we have a product of right-handed Dehn twists
 
 $$t_1\dots t_m=(t_1\dots t_{i_1})(t_{i_1+1}\dots t_2)\dots (t_{i_k+1}\dots t_m)\in \Gamma_1.$$
 
 This induces an elliptic fibration $X\rightarrow D^2$ having singular fibers $F_j$ ($j=1,\dots,k$), with monodromies conjugate to $(t_{i_j+1}\dots t_{i_{j+1}})$. If the product $t_1\dots t_m=1$ in $\Gamma_1$, then the induced fibration closes up to a fibration over $S^2$. Let $\Gamma_{1,k}$ denote the mapping class group of a compact, connected orientable genus $1$ surface with $k$ boundary components (the diffeomorphisms and isotopies are assumed to fix the boundary points). A lift of the relation $t_1\dots t_m=1$ to an identity in $\Gamma_{1,k}$ shows the existence of $k$ disjoint sections in the induced fibration. Let us denote by $t_{\delta_i}\in \Gamma_{1,k}$ the right-handed Dehn twist around a circle parallel to the $i$-th boundary component. If the lift to $\Gamma_{1,k}$ is of the form $t_1\dots t_m=t_{\delta_1}^{n_1}\dots t_{\delta_k}^{n_k}$ then the self-intersection of the $i$-th section is $-n_i$.

 In our case we will look at products of the form $t_1\dots t_{12}=1$ in $\Gamma_1$, which will give us an elliptic fibration $E(1)\rightarrow S^2$. Taking the $n$-th  power of the relation will give us an elliptic fibration $E(n)\rightarrow S^2$.
 
 \begin{remark}\label{rmk:monodromy}
    From this description, it is clear that given a factorisation of right-handed Dehn twists $t_1^{n_1}\dots t_m^{n_m}=1$ one can use the identity $t^{n+m}=t^n t^m$ to obtain another elliptic fibration which has one more discriminant point.
    Geometrically, this means for example that one $I_4$ fiber can be switched out for two $I_2$ fibers, without changing the diffeomorphism type of the total space, and the structure of the fibration is only changed near the original singular fiber.

    This shows that an $I_{n\times m}$ fiber can be replaced by $m$ number of $I_n$ fibers with the same monodromy. In particular, an $I_{4m}$ fiber can be replaced with $m$ number of $I_4$ fibers, and each $I_4$ can also be changed to $4$ nodal fibers, allowing two double node surgeries to be performed.
 \end{remark}

\subsection{An involution on $E(2n+1)$}\label{constr}
We briefly recall the construction of the extended involution, still denoted by $\iota:E(1)\to E(1)$ from \cite[Section 3]{stipsicz2023exotic}.
Consider $S^2\times S^2$ and the antipodal map on both factors, denoted by $p\times p$.
Let $a,c\in S^2$ and $b,d\in S^2$ such that $c\neq p(a)$ and $d\neq p(b)$. Let $C_0=S^2\times \{a\}\cup S^2\times \{c\}\cup \{b\}\times S^2\cup \{d\}\times S^2$.

This gives us four spheres (two fibers and two sections) in $S^2\times S^2$. Now $C_1=(p\times p)(C_0)$ is another configuration of four spheres with $|C_0\cap C_1|=8$.

\begin{figure}[ht!]
    \def\svgwidth{0,45\columnwidth}
    \centering
    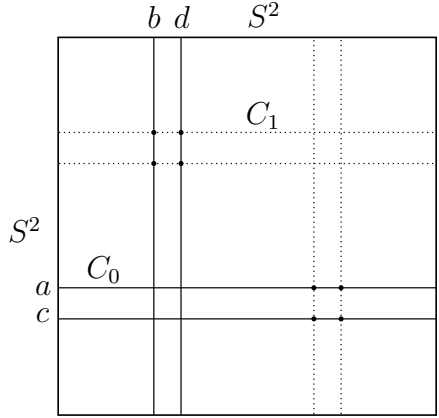

    \caption{$S^2\times S^2$ and the generators of the pencil $C_0, C_1$.}

\end{figure}

Note, that the involution $p\times p$ extends if we blow up $S^2\times S^2$ at pairs of these intersection points, thus getting a well defined fibration $S^2\times S^2\#8\overline{\CP^2}\to S^2$, still equipped with an involution, which we still call $\iota$.
This fibration is an elliptic fibration with total space $E(1)$.

We extend this involution further by taking a smooth fiber and its pair under $\iota$.
Taking a fiber sum with a copy of $E(n)$ along this fiber, and its pair, the involution extends by exchanging the two copies of $E(n)$ in the manifold $E(2n+1)=E(n)\#_f E(1)\#_f E(n)$.

\section{Construction}

Now putting the above results together, we will do a number of double node surgeries on some elliptic surface $E(2n+1)$ using the twist knots $K_m$ and $K_1$ depicted on Figure \ref{fig:twist}.
Here, $m$ denotes the parameter of our infinite exotic families on the fixed topological type, and $n$ is the coefficient in the main theorem which determines the topological type.
After the double node surgeries, we blow up the double points, and modify the fibration to produce a configuration $C_p$ which will be rationally blown down.

More precisely, consider the fibration $E(2n+1)$ described in Subsection \ref{constr}.
Choose a section $s\in\Gamma(E(2n+1))$, and consider the image $\iota(s)$ of that section under the involution map $\iota$ (they have self-intersection $-(2n+1)$).
Choose two smooth fibers (again equivariantly) of our fibration, and do knot surgery along them using the twist knot $K_m$.
\begin{figure}[ht!]
    \def\svgwidth{0,93\columnwidth}
    \centering
    %
\begingroup%
  \makeatletter%
  \providecommand\color[2][]{%
    \errmessage{(Inkscape) Color is used for the text in Inkscape, but the package 'color.sty' is not loaded}%
    \renewcommand\color[2][]{}%
  }%
  \providecommand\transparent[1]{%
    \errmessage{(Inkscape) Transparency is used (non-zero) for the text in Inkscape, but the package 'transparent.sty' is not loaded}%
    \renewcommand\transparent[1]{}%
  }%
  \providecommand\rotatebox[2]{#2}%
  \newcommand*\fsize{\dimexpr\f@size pt\relax}%
  \newcommand*\lineheight[1]{\fontsize{\fsize}{#1\fsize}\selectfont}%
  \ifx\svgwidth\undefined%
    \setlength{\unitlength}{742.65258606bp}%
    \ifx\svgscale\undefined%
      \relax%
    \else%
      \setlength{\unitlength}{\unitlength * \real{\svgscale}}%
    \fi%
  \else%
    \setlength{\unitlength}{\svgwidth}%
  \fi%
  \global\let\svgwidth\undefined%
  \global\let\svgscale\undefined%
  \makeatother%
  \begin{picture}(1,0.41629464)%
    \lineheight{1}%
    \setlength\tabcolsep{0pt}%
    \put(0,0){\includegraphics[width=\unitlength,page=1]{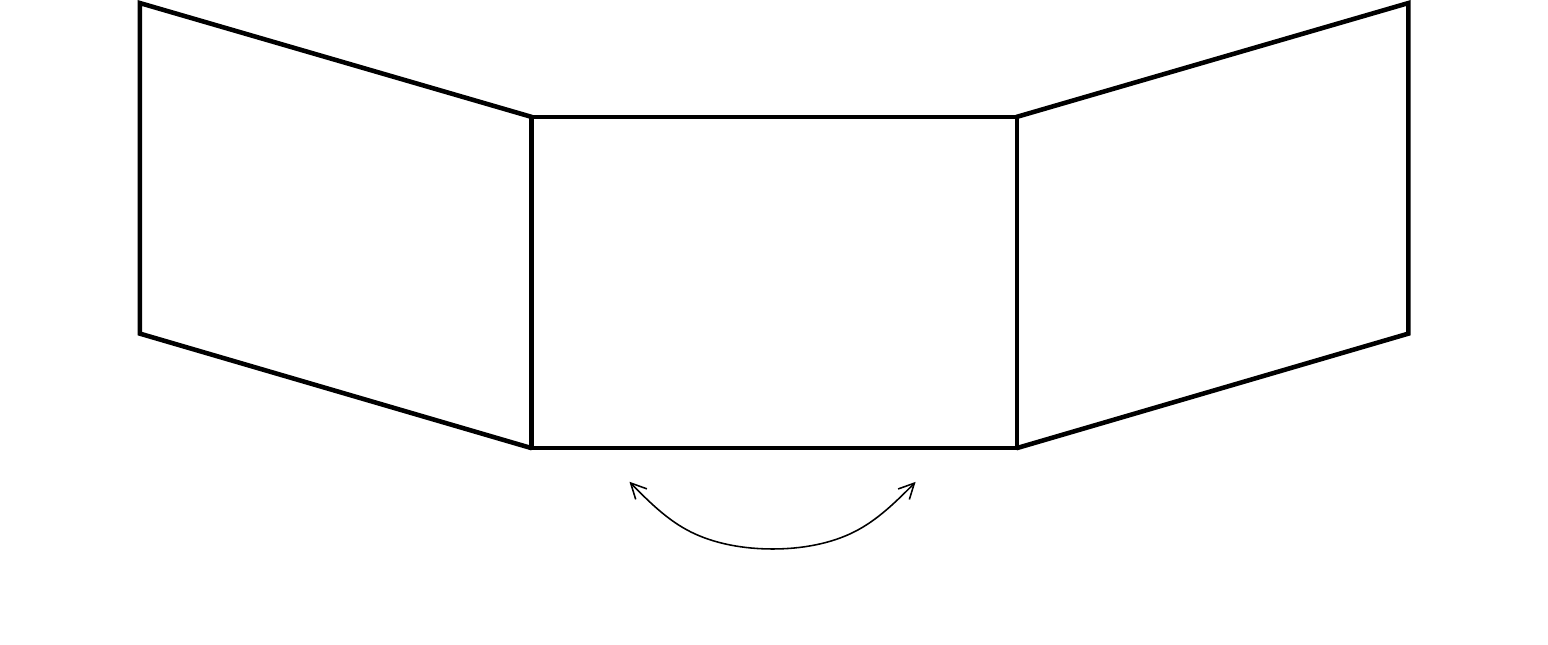}}%
    \put(0.43755403,0.03289979){\makebox(0,0)[lt]{\lineheight{1.25}\smash{\begin{tabular}[t]{l}$\mathbb{Z}/2\mathbb{Z}$-action\end{tabular}}}}%
    \put(0,0){\includegraphics[width=\unitlength,page=2]{elliptic.pdf}}%
    \put(0.92001476,0.26239796){\makebox(0,0)[lt]{\lineheight{1.25}\smash{\begin{tabular}[t]{l}$\iota (s)$\end{tabular}}}}%
    \put(0.92001476,0.34116971){\makebox(0,0)[lt]{\lineheight{1.25}\smash{\begin{tabular}[t]{l}$s$\end{tabular}}}}%
    \put(0,0){\includegraphics[width=\unitlength,page=3]{elliptic.pdf}}%
  \end{picture}%
\endgroup%

    \caption{Schematic picture of the elliptic surface $E(2n+1)$ seen as the fiber sum of $E(1)$ with two copies of $E(n)$.}
    \label{fig:elliptic}
\end{figure}

Since $g_3(K_m)=1$, the sections become genus two surfaces after the surgery, and following Section \ref{monodromy} we write
\begin{equation}\label{eq:nine}
1=(ab)^6=(a^3b)^3=a^4a^4a(b^{a^6}b^{a^3}b)
\end{equation}
in the monodromy group (see \cite[Section 2]{stipsicz2007singular} for further details) we can assume that there are two $I_4$ fibers in our $E(1)$, positioned symmetrically with respect to $\iota$.

\begin{remark}\label{rmk:conj}
    We use the relation $aba=bab$ (c.f. Subsection \ref{monodromy}) and then use conjugation $aba=a^2b^a$ several times to move the $a$'s to the required position. 
\end{remark}

\begin{figure}[ht!]
    \def\svgwidth{0,37\columnwidth}
    \centering
    %
\begingroup%
  \makeatletter%
  \providecommand\color[2][]{%
    \errmessage{(Inkscape) Color is used for the text in Inkscape, but the package 'color.sty' is not loaded}%
    \renewcommand\color[2][]{}%
  }%
  \providecommand\transparent[1]{%
    \errmessage{(Inkscape) Transparency is used (non-zero) for the text in Inkscape, but the package 'transparent.sty' is not loaded}%
    \renewcommand\transparent[1]{}%
  }%
  \providecommand\rotatebox[2]{#2}%
  \newcommand*\fsize{\dimexpr\f@size pt\relax}%
  \newcommand*\lineheight[1]{\fontsize{\fsize}{#1\fsize}\selectfont}%
  \ifx\svgwidth\undefined%
    \setlength{\unitlength}{374.54238177bp}%
    \ifx\svgscale\undefined%
      \relax%
    \else%
      \setlength{\unitlength}{\unitlength * \real{\svgscale}}%
    \fi%
  \else%
    \setlength{\unitlength}{\svgwidth}%
  \fi%
  \global\let\svgwidth\undefined%
  \global\let\svgscale\undefined%
  \makeatother%
  \begin{picture}(1,0.59408396)%
    \lineheight{1}%
    \setlength\tabcolsep{0pt}%
    \put(0,0){\includegraphics[width=\unitlength,page=1]{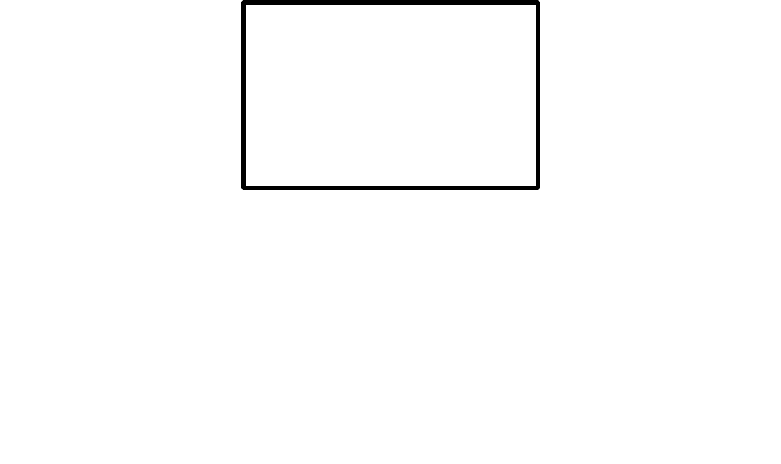}}%
    \put(0.38152033,0.44837917){\makebox(0,0)[lt]{\lineheight{1.25}\smash{\begin{tabular}[t]{l}$2m-1$\end{tabular}}}}%
    \put(0,0){\includegraphics[width=\unitlength,page=2]{twist-knot.pdf}}%
  \end{picture}%
\endgroup%

    \caption{The twist knot $K_m$, where the box represents ($2m-1$)-many half twists. The genus 1 Seifert surface is obtained from the immersed disc bounded by the two parallel strands.}
    \label{fig:twist}
\end{figure}
Using Remark \ref{rmk:monodromy} and \cite{fintushel2006double} we can perform double node surgeries on these $I_4$ fibers to 
exchange the genera for two positive double points on both sections.
Indeed, with one $I_4$ fiber positioned near the knot surgery we can trade one genus for a positive double point on both sections. The same can be done with two $I_4$ fibers equivariantly: carrying out the first two double node surgeries will change one genus on both $s$ and $\iota(s)$ to a positive double point, and doing the same construction on their "mirror image" changes the remaining genus into a positive double point on both sections (see Remark \ref{rmk:monodromy}).
Note that the blow up of a positive double point decreases the self-intersection of the blown up section by $4$, and decreases the homology class of the section by $2E$, where $E$ is the exceptional sphere of the blow up.
Blowing up the four double points, the sections now have self-intersection $-2n-9$.

Next, we do $4k$ additional double node surgeries along the trefoil knot (we choose the left-handed trefoil $K_1$) positioned symmetrically with respect to $\iota$ in the two $E(n)$ parts of our $E(2n+1)$.
\begin{remark}
Note that a number of other knots would work, but we choose the trefoil to keep the computations of the Seiberg-Witten invariants simple.
\end{remark}
As before, $g_3(K_1)=1$ and its Alexander polynomial is $\Delta_{K_1}=t-1+t^{-1}$, 
allowing us to apply additional double node surgeries to produce double points on the sections and blow them up.

\begin{figure}[ht!]
    \def\svgwidth{0,37\columnwidth}
    \centering
    %
\begingroup%
  \makeatletter%
  \providecommand\color[2][]{%
    \errmessage{(Inkscape) Color is used for the text in Inkscape, but the package 'color.sty' is not loaded}%
    \renewcommand\color[2][]{}%
  }%
  \providecommand\transparent[1]{%
    \errmessage{(Inkscape) Transparency is used (non-zero) for the text in Inkscape, but the package 'transparent.sty' is not loaded}%
    \renewcommand\transparent[1]{}%
  }%
  \providecommand\rotatebox[2]{#2}%
  \newcommand*\fsize{\dimexpr\f@size pt\relax}%
  \newcommand*\lineheight[1]{\fontsize{\fsize}{#1\fsize}\selectfont}%
  \ifx\svgwidth\undefined%
    \setlength{\unitlength}{482.67520183bp}%
    \ifx\svgscale\undefined%
      \relax%
    \else%
      \setlength{\unitlength}{\unitlength * \real{\svgscale}}%
    \fi%
  \else%
    \setlength{\unitlength}{\svgwidth}%
  \fi%
  \global\let\svgwidth\undefined%
  \global\let\svgscale\undefined%
  \makeatother%
  \begin{picture}(1,0.68494479)%
    \lineheight{1}%
    \setlength\tabcolsep{0pt}%
    \put(0,0){\includegraphics[width=\unitlength,page=1]{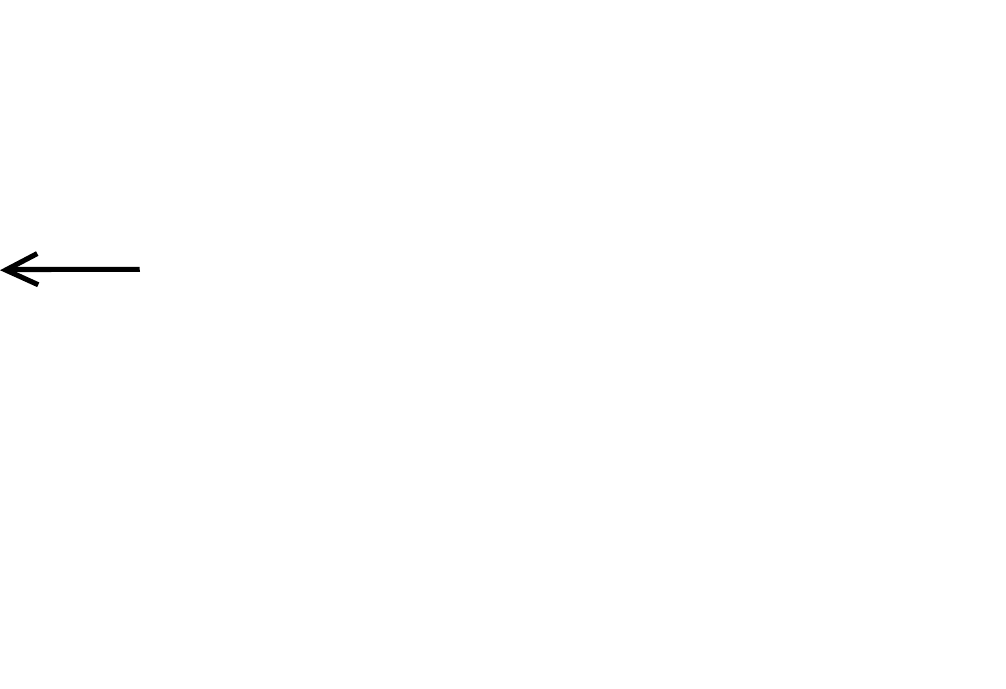}}%
    \put(0.043952,0.45095954){\makebox(0,0)[lt]{\lineheight{1.25}\smash{\begin{tabular}[t]{l}$s$\end{tabular}}}}%
    \put(0.763333,0.46037819){\makebox(0,0)[lt]{\lineheight{1.25}\smash{\begin{tabular}[t]{l}$\iota (s)$\end{tabular}}}}%
    \put(0,0){\includegraphics[width=\unitlength,page=2]{dual-graph.pdf}}%
    \put(0.37789557,0.00315624){\makebox(0,0)[lt]{\lineheight{1.25}\smash{\begin{tabular}[t]{l}$8n-3$\end{tabular}}}}%
    \put(0,0){\includegraphics[width=\unitlength,page=3]{dual-graph.pdf}}%
  \end{picture}%
\endgroup%

    \caption{Dual graph of the fiber with the two sections mapped to each-other.}
    \label{fig:dual}
\end{figure}

This lowers the self-intersection of both sections to $-2n-9-8k$.
Thus we need to find $k$ many $I_4$ fibers in $E(n)$ for this construction.
These fibers and their images under $\iota$ will be used to do the double node surgeries mentioned.
We produce a singular $I_{8n+1}$ fiber, where the dual 
graph of the fiber and the two sections looks like Figure \ref{fig:dual}.

\begin{figure}[ht!]
    \def\svgwidth{0,8\columnwidth}
    \centering
    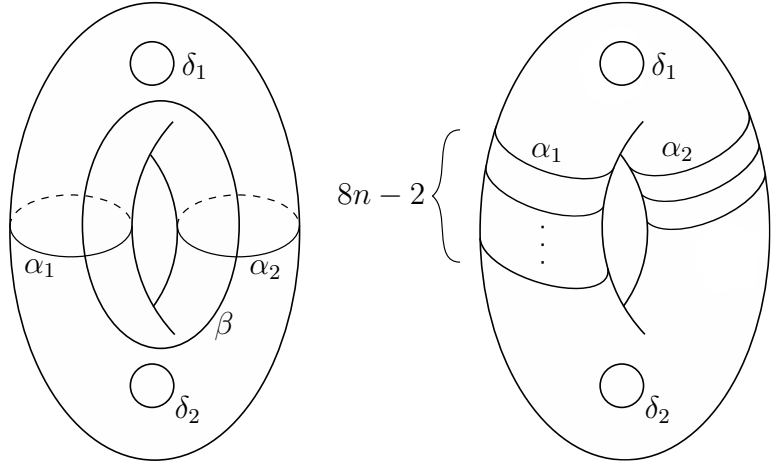

    \caption{Generators of the mapping class group of the twice punctured torus and the configuration produced by the factorisation.}
    \label{fig:two-holed-torus}
\end{figure}

The first equality of the next equation is a result of Korkmaz-Ozbagci (\cite[Section 3.2]{korkmaz2008sections}). From this, we compute the following factorisation in the mapping class group of the torus with two discs removed, where $\beta$ denotes a right handed Dehn twist along the longitude, and $\overline{\beta}$ is its inverse (see Figure \ref{fig:two-holed-torus}), using the same technique mentioned in Remark \ref{rmk:conj} to obtain the last two equalities:

\begin{align*}
\delta_1\delta_2=(\alpha_1\alpha_2\beta)^4=
(\alpha_1^3 \alpha_2^{\overline{\beta}}\alpha_1\beta)^2=\alpha_1^8\alpha_2^{\overline{\beta}\alpha_1^5}\beta^{\alpha_1^4}\alpha_2^{\overline{\beta}\alpha_1}\beta.
\end{align*}
\begin{remark} Note that Dehn twists along disjoint curves commute, hence $\alpha_1$ and $\alpha_2$ commute, and we can use cyclic permutation since $\alpha_1$, $\alpha_2$ and $\beta$ commute with $\delta_1$ and $\delta_2$. \end{remark}
As we can see, this decomposition does not contain any $\alpha_2$ (without conjugation). This means that in the fibration corresponding to this decomposition, the two sections would intersect the same $(-2)$-curve of the singular fiber.
In order to separate the two sections we also need the decomposition 
$\delta_1\delta_2=\alpha_1^6\alpha_2^3\beta^{\alpha_1^4\alpha_2^2}\beta^{\alpha_1^2\alpha_2}\beta$ which follows from the same principles as the previous.

Multiplying $n-1$ times the first decomposition with the second and using the method outlined in Section \ref{monodromy}, we obtain an identity of the form
\begin{align}\label{eq:factor}
\delta_1^n\delta_2^n=\alpha_1^{8n-2}\alpha_2^3\big((4n-1)\text{-many right handed Dehn twists}\big)
\end{align}
which implies that in the space $E(n)$ we can indeed achieve the configuration of Figure \ref{fig:dual}.

\begin{figure}[ht!]
    \def\svgwidth{1,2\columnwidth}
    \centering
    %
\begingroup%
  \makeatletter%
  \providecommand\color[2][]{%
    \errmessage{(Inkscape) Color is used for the text in Inkscape, but the package 'color.sty' is not loaded}%
    \renewcommand\color[2][]{}%
  }%
  \providecommand\transparent[1]{%
    \errmessage{(Inkscape) Transparency is used (non-zero) for the text in Inkscape, but the package 'transparent.sty' is not loaded}%
    \renewcommand\transparent[1]{}%
  }%
  \providecommand\rotatebox[2]{#2}%
  \newcommand*\fsize{\dimexpr\f@size pt\relax}%
  \newcommand*\lineheight[1]{\fontsize{\fsize}{#1\fsize}\selectfont}%
  \ifx\svgwidth\undefined%
    \setlength{\unitlength}{676.53826544bp}%
    \ifx\svgscale\undefined%
      \relax%
    \else%
      \setlength{\unitlength}{\unitlength * \real{\svgscale}}%
    \fi%
  \else%
    \setlength{\unitlength}{\svgwidth}%
  \fi%
  \global\let\svgwidth\undefined%
  \global\let\svgscale\undefined%
  \makeatother%
  \begin{picture}(1,0.39426795)%
    \lineheight{1}%
    \setlength\tabcolsep{0pt}%
    \put(0,0){\includegraphics[width=\unitlength,page=1]{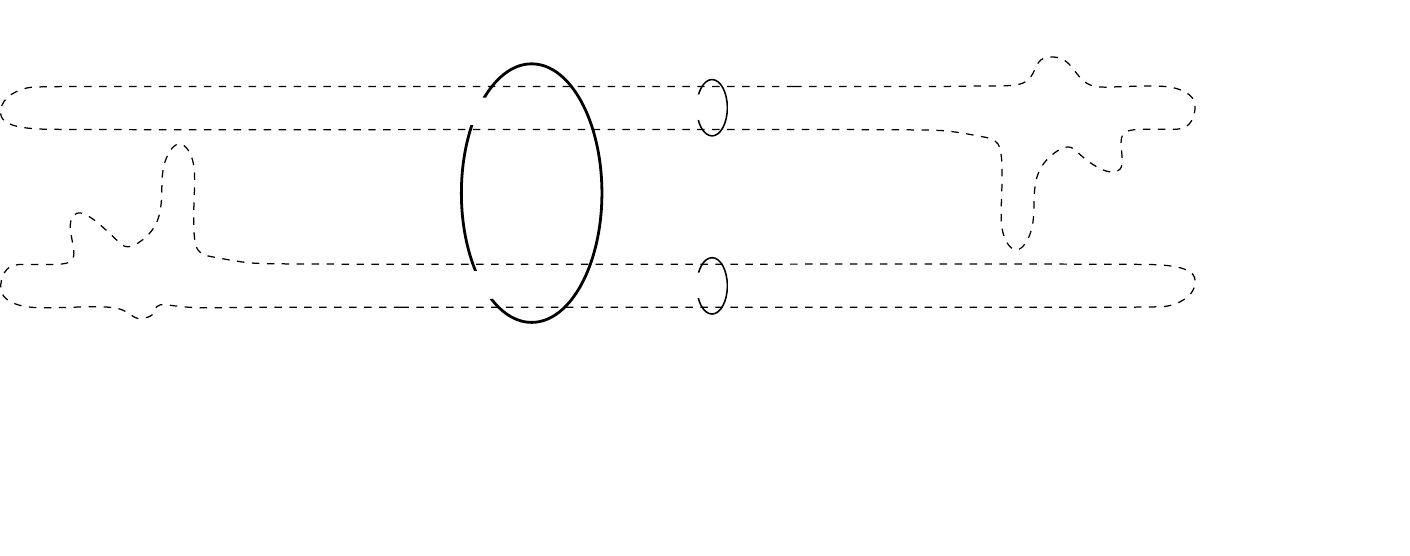}}%
    \put(0.85729548,0.31118795){\makebox(0,0)[lt]{\lineheight{1.25}\smash{\begin{tabular}[t]{l}$s$\end{tabular}}}}%
    \put(0.85691904,0.18788924){\makebox(0,0)[lt]{\lineheight{1.25}\smash{\begin{tabular}[t]{l}$\iota (s)$\end{tabular}}}}%
    \put(0.35260106,0.00505476){\makebox(0,0)[lt]{\lineheight{1.25}\smash{\begin{tabular}[t]{l}$\mathbb{Z}/2\mathbb{Z}$-action\end{tabular}}}}%
    \put(0,0){\includegraphics[width=\unitlength,page=2]{configuration.pdf}}%
    \put(0.16828794,0.24306775){\makebox(0,0)[lt]{\lineheight{1.25}\smash{\begin{tabular}[t]{l}$8n-3$\end{tabular}}}}%
    \put(0,0){\includegraphics[width=\unitlength,page=3]{configuration.pdf}}%
    \put(0.59841903,0.24306775){\makebox(0,0)[lt]{\lineheight{1.25}\smash{\begin{tabular}[t]{l}$8n-3$\end{tabular}}}}%
    \put(0,0){\includegraphics[width=\unitlength,page=4]{configuration.pdf}}%
  \end{picture}%
\endgroup%

    \caption{The discs in the complement of the configuration, and the punctured $S^2$'s (represented by thickened lines).}
    \label{fig:config-simply-conn}
\end{figure}

Consider the union of $s$, and the ($-2$)-curves of the just constructed {$I_{8n+1}$} and its pair under $\iota$ (see Figure \ref{fig:config-simply-conn}). Since $s$ has self-intersection $-(2n+8k+9)$, this contains a $C_{2n+7+8k}$ subconfiguration, which is disjoint from its pair under $\iota$ since we have two additional spheres in the fiber.
The fundamental group of the complement is generated by a normal circle of one of the endpoints of the linear plumbing $C_{2n+7+8k}$, the longer chain guarantees a bounding disc to this generator in the complement inside $E(2n+1)$.
Since $E(2n+1)$ is simply connected, any curve in $E(2n+1)\setminus (C_{2n+7+8k}\cup \iota(C_{2n+7+8k}))$ bounds a disc in the ambient space.
This intersects $\partial C_{2n+7+8k}=L((2n+7+8k)^2,2n+6+8k)$ in some arcs. In general, lens spaces have abelian fundamental group, so it is enough to calculate homologically. These intersection arcs can be homotoped to the end of the chain, where they bound a disc.

Consider one of the $E(n)$ submanifolds and its pair under $\iota$. We wish to compute a bound on $k$ and to this end give an overview of our construction.
It is enough to make the computation in one of the $E(n)$ parts of $E(2n+1)$, since its pair under $\iota$ is disjoint from it, we can easily extend everything equivariantly.
To increase $k$ by one, we need an $I_4$ fiber to do two double node surgeries, and after blowing up each section once, their framing drops by $8$, forcing us to construct a $C_{p+8}$ instead of a $C_p$ (the $C_p$ where the first vertex is $s$ has its $-2$ curves in $E(n)$, and its pair has its $-2$ curves in $\iota(E(n))$).
Thus all together increasing $k$ once requires $\alpha_1^{12}$ and when $k=0$ we need an $\alpha_1^{2n+7}$, since the double node surgeries were done in the "middle" $E(1)$ part, where the involution acts non-trivially.
By Equation (\ref{eq:factor}), we can produce $\alpha_1^{8n-2}$ in the monodromy factorization of $E(n)$. 
Therefore, we have 
\begin{equation}\label{eq:knumber}
8n-2=2n+7+8k+4k,
\end{equation}
i.e. $k=\lfloor\frac n2-\frac3{4}\rfloor$ gives the maximum number of times we can perform the surgery. Note that this means the coefficient $8n-6k$ is equal to $5n+6$ if $n$ is even, and to $5n+9$ if $n$ is odd.

\begin{remark}
    {In the process, each double node surgery needs a pair of Lefschetz thimbles (corresponding to parallel vanishing cycles), and different surgeries need separate thimbles to achieve that the end-result (after the blow-ups) is embedded, rather than immersed. This explains the need of the large number of parallel vanishing cycles in our construction, and in turn the coefficient of $k$ in Equation (\ref{eq:knumber}).}

\end{remark}

Now, we are ready to prove the main result:
\begin{theorem}
\label{maintheorem}
    There are infinitely many {irreducible} smooth structures on $Z_1\#$ $2n\CP^2\#$ $(8n-6k)\overline{\CP^2}$ for $n\in\mathbb{N}$ and $0\leq k\leq\lfloor \frac n2-\frac{3}{4}\rfloor$ where $4\not|\ n-k$.
\end{theorem}
\begin{remark}
In order to show that our manifolds are non-spin, we use Rohlin's theorem (a spin four-manifold has signature divisible by $16$, \cite{rohlin1952new}). 
By the calculation below this means that we have to exclude those cases where $n\equiv 5$ or $6 \pmod{8}$ (\ie when $n-k$ is divisible by 4).
Note that the construction remains valid in these cases, and we suspect that the manifolds are non-spin as well.
\end{remark}
\begin{proof}

The Seiberg-Witten invariant of $E(2n+1)$ is
\[
\mathcal{SW}_{E(2n+1)}=(e^{f}-e^{-f})^{2n-1}
\]
where $f\in H^2(E(2n+1);\Z)$ is the fiber class (\cite{fintushel1997rational} \cite[Example 1]{akhmedov2007construction}).
By \cite{fintushel1998knots}, we get that after both surgeries with $K_m$ the invariants get multiplied by $me^{2f}-(2m-1)+me^{-2f}$, and by $e^{2f}-1+e^{-2f}$ after each additional $K_1$-surgery.
Since $E(2n+1)$ is of simple type, the potential basic classes become $(2n-1-2r)f$, where $r\in\{-2k-2,\dots,2n+2k+1\}$ (and this new manifold, homeomorphic to $E(2n+1)$ is still of simple type).

After the $4k+4$ blow ups the basic classes are of the form $(2n-1-2r)f+\sum_{i=1}^{4k+4}\pm E_i$ with $r$ as before.
The sections get blown up at double points, so the sphere is represented by $[s]-\sum_{i=1}^{2k+2}2E_i$ and its pair by $[\iota(s)]-\sum_{i=2k+3}^{4k+4}2E_i$. Note that the construction is done in an equivariant manner on two disjoint sections.

We compute:
\[\big\langle (2n-1-2r)f+\sum_{i=1}^{4k+4}\pm E_i,[s]-\sum_{i=1}^{2k+2}2E_i\big\rangle=2n-1-2r+2a,\]
where $a$ represents the difference of the number of $E_i$'s with negative sign in the basic class and the number of $E_i$'s with positive sign for $1\leq i\leq 2k+2$, so it is any even number satisfying $|a|\leq 2k+2$ and $\langle\cdot,\cdot\rangle$ is the intersection pairing.
Thus the value above is at least $2n-1-4n-2-4k-4k-4=-2n-8k-7$ and at most $2n-1+4k+4+4k+4=2n+8k+7$ meaning that both configurations are tautly embedded (see Definition \ref{def:taut}), since the ($-2$)-spheres of the configuration are in a fiber, so $f$ evaluates on them as $0$.
 
We apply \cite[Theorem 8.5.18.]{gompf19994} to get that the only basic classes which extend to the rational blow down are $\alpha=(2n-1+4k+4)f+\sum_{i=1}^{4k+4}E_i$ and its negative, since this is the only class which evaluates maximally on both sections.
This class corresponds to the leading coefficient of the invariant before the rational blow down.
The leading coefficient of the product $$(t-t^{-1})^{2n-1}(mt^2-(2m-1)+mt^{-2})^2(t^2-1+t^{-2})^{2k}$$ is $m^2$, thus these manifolds are all smoothly distinct.

Topologically, $\chi(E(2n+1))=24n+12$, which we change by $4k+4-4n-12-16k$ with the blow ups and the rational blow downs to obtain $20n-12k+4$.
Furthermore, $b_2^+(E(2n+1))=4n+1,\ b_2^-(E(2n+1))=20n+9$, we only add and remove negative definite submanifolds, so the signature becomes $4n+1-(20n+9+4k+4-4n-12-16k)=-12n+12k$.
This number is not divisible by 16 by assumption, thus our manifolds are not spin \ie type I.
Taking the quotient by $\iota$ halves both $\sigma$ and $\chi$, so by the homeomorphism classification of $\pi_1=\mathbb{Z}/2\mathbb{Z}$ smooth manifolds the factor is homeomorphic to $Z_1\#2n\CP^2\#(8n-6k)\overline{\CP^2}$, but as detected by the Seiberg-Witten invariant of its universal cover, they are not diffeomorphic.

Finally by the above calculation and Equation~\eqref{eq:smiple}, we see that $\alpha^2={3\sigma+2\chi}=4n+12k+8$ for our basic class $\alpha$.
Irreducibility follows from this, since in a reducible manifold there are two basic classes, the difference of which has square $-4$ by \cite[Lemma 2.3]{stipsicz2023definite}, but in our case this number is always positive.
This shows that the universal cover of our manifolds are irreducible, and thus the quotient is also irreducible.
\end{proof}

\nocite{*}
\bibliographystyle{plain}
\bibliography{references}

\bigskip
\footnotesize

M\'{a}rton Beke$^{1,3}$
\textit{E-mail address}: \texttt{bekem@renyi.hu}

\medskip

L\'{a}szl\'{o} Koltai$^{2,3}$
\textit{E-mail address}: \texttt{koltail@renyi.hu}

\medskip

Sarah~Zampa$^{1,3}$
\textit{E-mail address}: \texttt{zampa.sarah@renyi.hu}

\bigskip
\textsc{$^{1}$ Department of Geometry and Algebra, Institute of Mathematics, Budapest University of Technology and Economics, Műegyetem rkp. 3., H-1111 Budapest, Hungary}\par
\textsc{$^{2}$ Department of Analysis, Institute of Mathematics, E\"{o}tv\"{o}s Lor\'{a}nd University, Pázmány Péter sétány 1/C., H-1117 Budapest, Hungary}\par
\textsc{$^{3}$ Department of Algebraic Geometry and Differential Topology, HUN-REN Alfréd R\'{e}nyi Institute of Mathematics, Reáltanoda u. 13-15., H-1053 Budapest, Hungary}

\end{document}